\documentclass[twoside,draft,reqno]{amsart}
\usepackage{amssymb}
\numberwithin{equation}{section}

\theoremstyle{plain}
\newtheorem{Thm}{Theorem}[section]
\newtheorem{Cor}[Thm]{Corollary}
\newtheorem{Lem}[Thm]{Lemma}
\newtheorem{Prop}[Thm]{Proposition}

\theoremstyle{definition}
\newtheorem{Def}[Thm]{Definition}

\theoremstyle{remark}

\newtheorem{Rem}[Thm]{Remark}

\newif\ifShowLabels
\ShowLabelstrue
\newdimen\theight
\def\TeXref#1{%
        \leavevmode\vadjust{\setbox0=\hbox{{\tt
                \quad\quad  {\small \textrm #1}}}%
        \theight=\ht0
        \advance\theight by \lineskip
        \kern -\theight \vbox to
        \theight{\rightline{\rlap{\box0}}%
        \vss}%
        }}%

\ShowLabelsfalse

\renewcommand{\sec}[2]{\section{#2}\label{S:#1}%
        \ifShowLabels \TeXref{{S:#1}} \fi}
\newcommand{\ssec}[2]{\subsection{#2}\label{SS:#1}%
        \ifShowLabels \TeXref{{SS:#1}} \fi}

\newcommand{\refs}[1]{Section ~\ref{S:#1}}
\newcommand{\refss}[1]{Subsection ~\ref{SS:#1}}
\newcommand{\reft}[1]{Theorem ~\ref{T:#1}}
\newcommand{\refl}[1]{Lemma ~\ref{L:#1}}

\newcommand{\refc}[1]{Corollary ~\ref{C:#1}}

\newcommand{\refe}[1]{\eqref{E:#1}}

\newenvironment{thm}[1]%
        { \begin{Thm} \label{T:#1}  \ifShowLabels \TeXref{T:#1} \fi }%
        { \end{Thm} }

\renewcommand{\th}[1]{\begin{thm}{#1} \sl }
\renewcommand{\eth}{\end{thm} }

\newenvironment{lemma}[1]%
        { \begin{Lem} \label{L:#1}  \ifShowLabels \TeXref{L:#1} \fi }%
        { \end{Lem} }
\newcommand{\lem}[1]{\begin{lemma}{#1} \sl}
\newcommand{\elem}{\end{lemma}}

\newenvironment{propos}[1]%
        { \begin{Prop} \label{P:#1}  \ifShowLabels \TeXref{P:#1} \fi }%
        { \end{Prop} }
\newcommand{\prop}[1]{\begin{propos}{#1}\sl }
\newcommand{\eprop}{\end{propos}}

\newenvironment{corol}[1]%
        { \begin{Cor} \label{C:#1}  \ifShowLabels \TeXref{C:#1} \fi }%
        { \end{Cor} }
\newcommand{\cor}[1]{\begin{corol}{#1} \sl }
\newcommand{\ecor}{\end{corol}}

\newenvironment{defeni}[1]%
        { \begin{Def} \label{D:#1}  \ifShowLabels \TeXref{D:#1} \fi }%
        { \end{Def} }
\newcommand{\defe}[1]{\begin{defeni}{#1} \sl }
\newcommand{\edefe}{\end{defeni}}

\newenvironment{remark}[1]%
        { \begin{Rem} \label{R:#1}  \ifShowLabels \TeXref{R:#1} \fi }%
        { \end{Rem} }
\newcommand{\rem}[1]{\begin{remark}{#1}}
\newcommand{\erem}{\end{remark}}

\newcommand{\eq}[1]%
        { \ifShowLabels \TeXref{E:#1} \fi
           \begin{equation} \label{E:#1} }
\newcommand{\eeq}{\end{equation}}

\newcommand{\prf}{ \begin{proof} }
\newcommand{\eprf}{ \end{proof} }

\newcommand\alp{\alpha}         
\newcommand\bet{\beta}
\newcommand\gam{\gamma}         
         
\newcommand\eps{\varepsilon}

\newcommand\iot{\iota}

\newcommand\lam{\lambda}                
\newcommand\sig{\sigma}         

\newcommand\ome{\omega}         \newcommand\Ome{\Omega}

\newcommand\calH{{\mathcal{H}}}
\newcommand\calL{{\mathcal{L}}}
\newcommand\calO{{\mathcal{O}}}
\newcommand\calU{{\mathcal{U}}}

\newcommand\RR{\mathbb{R}}
\newcommand\CC{\mathbb{C}}
\newcommand\NN{\mathbb{N}}

\newcommand\grk{{\mathfrak{k}}}

\newcommand\nek{,\ldots,}
\newcommand\sdp{\times \hskip -0.3em {\raise 0.3ex
\hbox{$\scriptscriptstyle |$}}} 

\newcommand\Ad{\operatorname{Ad\, }}
\newcommand\diag{\operatorname{diag}}

\newcommand\Id{\operatorname {Id}}
\newcommand\Ker{\operatorname{Ker}}
\newcommand\supp{\operatorname{supp}}

\newcommand\Vol{\operatorname{Vol}}

\newcommand\oE{{\overline{E}}}
\newcommand\oq{{\overline{q}}}

\newcommand\tilh{{\widetilde{h}}}

\renewcommand{\>}{\rangle}
\newcommand{\<}{\langle}

\renewcommand{\d}{\text{\( \partial\)}}\newcommand{\p}{\bar{\d}}

\newcommand{\w}{\omega}
\newcommand{\Xss}{X^{ss}}\newcommand{\Lss}{L^{ss}}
\newcommand{\X}{X_0}\renewcommand{\L}{L_0}

\newcommand{\n}{\nabla}
\newcommand{\hm}[1]{{H^{#1}(X,\calO(L))}}
\newcommand{\hmss}[1]{{H^{#1}(\Xss,\calO(\Lss))}}
\newcommand{\hmo}[1]{{H^{#1}(\X,\calO(\L))}}

\newcommand{\lap}{\bar\square_t}

\newcommand{\h}[1]{\calH^{#1}(X_0,L_0)}

\newcommand{\ka}{K\"ahler }
\begin{document}
\hskip5.5cm{To appear in ``Quantization of Singular Symplectic
Quotients"}

\

\title{Cohomology of the Mumford Quotient}
\author{Maxim Braverman}
\address{Department of Mathematics\\
        Northeastern University   \\
        Boston, MA 02115 \\
        USA
         }
\email{maxim@neu.edu}
\thanks{This research was partially supported by grant No. 96-00210/1 from
the United States-Israel Binational Science Foundation (BSF)}

\begin{abstract}
Let $X$ be a smooth projective variety acted on by a reductive
group $G$. Let $L$ be a positive $G$-equivariant line bundle over
$X$. We use a Witten type deformation of the Dolbeault complex of
$L$, introduced by Tian and Zhang, to show, that the cohomology of
the sheaf of holomorphic sections of the induced bundle on the
Mumford quotient of $(X,L)$ is equal to the $G$-invariant part on
the cohomology of the sheaf of holomorphic sections of $L$. This
result, which was recently proven by C.~Teleman by a completely
different method, generalizes a theorem of Guillemin and
Sternberg, which addressed the global sections. It also shows,
that the Morse-type inequalities of Tian and Zhang \cite{TianZh98}
for symplectic reduction are, in fact, equalities.
\end{abstract}
\maketitle

\sec{introd}{Introduction}

\ssec{ss}{Cohomology of the set of semistable points}
Suppose $X$ is a smooth complex projective variety endowed with a
holomorphic action of a complex reductive group $G$. Given a
positive $G$-equivariant holomorphic line bundle $L$ over $X$,
Mumford \cite{MFK94} defined the notion of stability: A point
$x\in X$ is called {\em semi-stable} \/ with respect to $L$ if and
only if there exist a positive integer $m\in \NN$ and an invariant
section $s:X\to L^m$ such that $s(x)\not=0$. Let $\Xss$ denote the
set of semi-stable points of $X$. It is a $G$-invariant Zariski
open subset of $X$. Let $\Lss$ denote the restriction of $L$ to
$\Xss$.

The group $G$ acts naturally on the cohomology $\hm{*}$ and
$\hmss{*}$ of the sheaves of algebraic sections of $L$ and $\Lss$
respectively. We denote by $\hm{*}^G$ and $\hmss{*}^G$ the
subspaces of $G$-invariant elements in $\hm{*}$ and $\hmss{*}$.

The main result of this paper is the following
\th{mainA}
Assume that the action of $G$ on $\Xss$ is free. Then
\begin{equation}\label{E:mainA}
        \dim_\CC \hm{j}^G \ = \ \dim_\CC \hmss{j}^G,
                        \qquad j=0,1,\dots
\end{equation}
\eth
\rem{isomorph}
The equality \refe{mainA}, implies that there is an isomorphism
between vector spaces $\hm{j}^G$ and $\hmss{j}^G$. We actually
construct a family of such isomorphisms in the course of the
proof. Unfortunately, we don't have a way to single out one
canonical isomorphism.
\erem

\subsection{Historical remarks}\label{SS:hist}
For $j=0$, the equality \refe{mainA} was established by
Guillemin and Sternberg \cite[\S5]{GuiSter82}. Guillemin and
Sternberg also conjectured that
\eq{GS}
    \sum_{j=0}^n\, (-1)^j\dim_\CC\hm{j}^G \ = \
    \sum_{j=0}^n\, (-1)^j \dim_\CC\hmss{j}^G.
\end{equation}
This was proven by Meinrenken in \cite{Meinr-GS}. An analytic
proof of \refe{GS} was given by Tian and Zhang
\cite[Th.~0.4]{TianZh98}, who also proved the inequality
\eq{ineq}
    \dim_\CC\hm{j}^G\le \dim_\CC\hmss{j}^G.
\end{equation}
\reft{mainA} may be considered as an amplification of the results
of Guillemin-Sternberg, Meinrenken and Tian-Zhang.

The equality \refe{mainA} (with somehow weaker assumptions on the
action of $G$ on $X$) was recently proven by C.~Teleman
\cite{Telem98}. His proof is completely algebraic. It is based on
the study of the sheaf cohomology of $\calO(L)$ with support on a
stratum of the Morse stratification for the square of the moment
map (see also \cite{Ramadas96}, where a similar technique was
used).

The purpose of this paper is to give a direct analytic proof of
\reft{mainC}, based on the study of the Witten type deformation of
the Dolbeault complex of $X$ with values in $L$,  introduced by
Tian and Zhang  \cite{TianZh98}.

Certain generalizations of  \refe{GS}  for ``singular reduction",
i.e., when the action of $G$ on $\Xss$ is not free, was obtained
by Sjamaar \cite{Sj95}, Meinrenken and Sjamaar \cite{MeinSj99}.
Tian and Zhang \cite{TianZh98} gave a new analytic proof of these
results and also extended the inequalities \refe{ineq} to the
singular reduction. Recently, Zhang \cite{Zhang-sing} used the
method of this paper to extend \refe{mainA} to singular reduction.

\ssec{Mquot}{The cohomology of the Mumford quotient}
Suppose again that the action of $G$ on the set $\Xss$ of
semi-stable points is free. Then the quotient space $\Xss/G$ is a
smooth projective variety, called the {\em Mumford quotient} \/ of
$X$. The quotient $\Lss/G$ has a natural structure of a
holomorphic line bundle over $\Xss/G$. Clearly, the quotient map
$q:\Xss\to\Xss/G$ induces a natural isomorphism
\begin{equation}\label{E:Xs=X0}
        H^j(\Xss/G,\calO(\Lss/G)) \ \simeq \hmss{j}^G, \quad
                                        j=0,1,\dots
\end{equation}
Hence, \reft{mainA} is equivalent to the following
\th{mainB}
Assume that the action of $G$ on $\Xss$ is free. Then
\begin{equation}\label{E:mainB}
        \dim_\CC \hm{j}^G \ = \ \dim_\CC H^j(\Xss/G,\calO(\Lss/G)),
                        \qquad j=0,1,\dots
\end{equation}
\eth

\ssec{sympred}{Reformulation in terms of symplectic reduction}
Let $K\subset G$ be a maximal compact subgroup of $G$ and let
$\grk=Lie(K)$ denote the Lie algebra of $K$. Since, $G$ is a
complexification of $K$, it is clear, that the space of
$G$-invariant elements of $\hm{*}$ coincides with the space of
$K$-invariant elements:
\begin{equation}\label{E:G=K}
        \hm{j}^G \ = \ \hm{j}^K, \qquad j=0,1,\dots
\end{equation}

Fix a $K$-invariant hermitian structure on $L$ and let $\n$ be a
$K$-invariant holomorphic hermitian connection on $L$. Denote
\[
        \w \ = \ \frac{i}{2\pi}\, (\n)^2.
\]
Then $\w$ is a $K$-invariant \ka form on $X$, representing the
Chern class of $L$ in the integer cohomology of $X$.

For any section $s:X\to L$ and any $V\in \grk$, we denote by
$\calL_V s$ the infinitesimal action of $V$ on $s$, induced by the
action of $K\subset G$ on $L$.

Define a map $\mu:X\to \grk^*$ by the formula
\[
        \<\, \mu,\, V\, \> \ = \ \frac{i}{2\pi}\, (\calL_V \, - \, \n_V),
        \qquad  V\in \grk.
\]
Here $\n_V$ denotes the covariant derivative along the vector
field generated by $V$ on $X$. Then, cf. \cite{Kostant70}, $\mu$
is a {\em moment map} \/ for the action of $K$ on the symplectic
manifold $(X,\w)$.

The assumption that $G$ acts freely on $\Xss$ is equivalent to the
statement that $0\in \grk^*$ is a regular value of $\mu$ and $K$
acts freely on $\mu^{-1}(0)$, cf. \cite[\S4]{GuiSter82}. Let
$\X=\mu^{-1}(0)/K$ denote the symplectic reduction of $X$ at $0$.
Then, cf. \cite[\S5]{GuiSter82}, $\X$ is complex isomorphic to the
Mumford quotient:
\begin{equation}\label{E:X0=mum}
        \X \ = \ \Xss/G.
\end{equation}
Let $\oq:\mu^{-1}(0)\to \mu^{-1}(0)/K$ denote the quotient map.
Then, cf. \cite[Th.~3.2]{GuiSter82}, there exists a unique
holomorphic line bundle $\L$ on $\X$ such that $\oq^*\L=
L|_{\mu^{-1}(0)}$. Moreover, under the isomorphism \refe{X0=mum},
$\L$ is isomorphic to $\Lss/G$. Hence, it follows from \refe{G=K},
that Theorems~\ref{T:mainA}, \ref{T:mainB} are equivalent to the
following
\th{mainC}
Suppose $0$ is a regular value of $\mu$ and $K$ acts freely on
$\mu^{-1}(0)$. Then
\begin{equation}\label{E:mainC}
        \dim_\CC \hm{j}^K \ = \ \dim_\CC \hmo{j},
                        \qquad j=0,1,\dots
\end{equation}
\eth
The proof of \reft{mainC} is given in \refss{prmainC}. Here we
will only explain the main idea of the proof.

\ssec{sketch}{A sketch of the proof}
Following Tian and Zhang \cite{TianZh98}, we consider the one
parameter family \linebreak$\p_t=e^{-t\frac{|\mu|^2}2}\p
e^{t\frac{|\mu|^2}2}$ of differentials on the Dolbeault complex
$\Ome^{0,*}(X,L)$. Let $\lap=\p_t\p_t^*+\p_t^*\p_t$ denote the
corresponding family of Laplacians and let $\lap^{j,K}$ denote the
restriction of $\lap$ to the space of $K$-invariant forms in
$\Ome^{0,j}(X,L)$. Let $E^{j,K}_{t,\lam}$ denote the span of
eigenforms of $\lap^{j,K}$ with eigenvalues smaller than $\lam$.
Then $H^*(X,\calO(L))^K$ is isomorphic to the cohomology of the
complex $(E^{*,K}_{t,\lam},\p_t)$.

In \cite{TianZh98}, Tian and Zhang showed that there exists
$\lam>0$ such that
\[
    \dim H^j(X_0,\calO(L_0))=\dim E^{j,K}_{t,\lam},
        \qquad j=0, 1,\dots
\]
for any  $t\gg0$. Moreover, if we denote by
$\calH^j(X_0,L_0)\subset\Ome^{0,j}(X_0,L_0)$ the subset of
harmonic forms (with respect to the metrics introduced in
\refss{pr}), Tian and Zhang constructed an explicit isomorphism of
vector spaces
\[
    \Phi^j_{\lam,t}:\calH^j(X_0,L_0)\to E^{j,K}_{t,\lam}.
\]

Recall that we denote by $q$ the quotient map $\Xss\to
X_0=\Xss/G$. In \refs{integr}, we consider the integration map
$I_t:\Ome^{0,j}(X,L)\to\Ome^{0,j}(X_0,L_0)$ defined by the formula
\footnote{The map $I$ is defined as an integral over the fibers of $q$.
One can think about $I$ as a kind of equivariant push-forward of
differential forms under $q$.}
\[
        I_t:\, \alp \ \mapsto \
            \Big(\frac{t}{2\pi}\Big)^{r/4}\, \int_{q^{-1}(x)}\,
                e^{-\frac{t|\mu|^2}2}\alp\wedge\ome^r,
\]
where $r=\dim_\CC G$. Then we show (cf. \reft{integr}) that
$I_t\p_t=\p I_t$, for any $t\ge0$. Also the restriction of $I_t$
to $E_{\lam,t}^{*,K}$ is ``almost equal" to
$(\Phi^*_{\lam,t})^{-1}$ for $t\gg0$.

It follows from the existence of the map $I_t$ with the above
properties, that the restriction of $\p_t$ to $E_{\lam,t}^{*,K}$
vanishes (cf. \refc{integr}). This implies \reft{mainC} and,
hence, Theorems~\ref{T:mainA} and \ref{T:mainB}.

\subsection*{Acknowledgments}
I would like to thank Yael Karshon for valuable and stimulating
discussions. I am also very grateful to Weiping Zhang and the
referees for careful reading of the manuscript and  suggesting a
number of important comments and corrections.

\sec{deform}{Witten type deformation of the Dolbeault complex. The Tian-Zhang theorem}

Our proof of \reft{mainC} uses the Witten type deformation of the
Dolbeault complex, introduced by Tian and Zhang \cite{TianZh98}.
For convenience of the reader, we review in this section the
results of Tian and Zhang which will be used in the subsequent
sections.

\ssec{dt}{Deformation of the Doulbeault complex}
Let $\Ome^{0,j}(X, L)$ denote the space of smooth
$(0,j)$-differential forms on $X$ with values in $L$. The group
$K$ acts on this space and we denote by $\Ome^{0,j}(X, L)^K$ the
space of $K$-invariant elements in $\Ome^{0,j}(X, L)$.

Recall that $\mu:X\to\grk^*$ denotes the moment map for the $K$
action on $X$. Let $\grk$ (and, hence, $\grk^*$) be equipped with
an $\Ad K$-invariant metric, such that the volume of $K$ with
respect to the induced Riemannian metric is equal to 1. Let
$h_1\nek h_r$ denote an orthonormal basis of $\grk^*$. Denote by
$v_i$ the Killing vector fields on $X$ induced by the duals of
$h_i$. Then $\mu=\sum\,{}\mu_ih_i$, where each $\mu_i$ is a real
valued function on $X$ such that $\iot_{v_i}\ome=d\mu_i$. Hence,
\begin{equation}\label{E:dmu^2}
        d\, \frac{|\mu|^2}2 \ = \ \sum_{i=1}^r\, \mu_i d\mu_i
                        \ = \  \sum_{i=1}^r\, \mu_i \iot_{v_i}\ome.
\end{equation}

Consider the one-parameter family $\p_t, \ t\in\RR$ of
differentials on $\Ome^{0,j}(X, L)$, defined by the formula
\[
        \p_t:\, \alp \ \mapsto e^{-t\frac{|\mu|^2}2}\,
                         \p\, e^{t\frac{|\mu|^2}2}\, \alp
                \ = \ \p\alp \ + \ t\, \sum_{i=1}^k\, \mu_i \p\mu_i\wedge\alp.
\]
Clearly, $\p_t^2=0$ and the cohomology of the complex
$(\Ome^{0,*}(X, L),\p_t)$ is $K$-equivari\-ant\-ly isomorphic to
$H^*(X,\calO(L))$, for any $t\in\RR$.

\ssec{lapt}{The deformed Laplacian}
Fix a $K$-invariant \ka metric on $X$ and a $K$-invariant
Hermitian metric on $L$. Let $\p_t^*$ denote the formal adjoint of
$\p_t$ with respect to these metrics and set
\[
        \lap \ = \ \p^*_t\p_t \ + \ \p_t \p_t^*.
\]
Let $\lap^{j,K}$ denote the restriction of $\lap$ to the space
$\Ome^{0,j}(X, L)^K$.

For any $j=0,1,\dots, \lam>0$ and $t>0$, we define
$E_{\lam,t}^{j,K}$ to be the span of the eigenforms of
$\lap^{j,K}$ with eigenvalues less or equal than $\lam$. Then
$E^{*,K}_{\lam,t}$ is a subcomplex of $(\Ome^{0,*}(X, L),\p_t)$.
\, Since $E^{*,K}_{\lam,t}$ contains the kernel of $\lap$, it
follows from the Hodge theory that the cohomology of
$E^{*,K}_{\lam,t}$ is isomorphic to $H^*(X,\calO(L))^K$.

The following theorem of Tian and Zhang
\cite{TianZh98}\footnote{\reft{TZ1} is a combination of Th.~3.13
and \S4.d of \cite{TianZh98}} is crucial for our paper:
\th{TZ1}
There exist $\lam, t_0>0$ such that, for any $j=0,1,\ldots$ and
any $t>t_0$, \, $\lam$ is not in the spectrum of \/ $\lap^{j,K}$
and
\begin{equation}\label{E:TZ1}
        \dim E^{j,K}_{\lam,t} \ = \ \dim H^j(X_0,\calO(L_0)).
\end{equation}
\eth
As an immediate consequence of \reft{TZ1} and the fact that
$H^*(X,\calO(L))^K$ is isomorphic to the cohomology of the complex
$(E^{*,K}_{\lam,t},\p_t)$ we obtain the following inequalities of
Tian and Zhang \cite[Th.~4.8]{TianZh98}:
\begin{equation}\label{E:Tzin}
        \dim H^j(X,\calO(L))^K \ \le \ \dim H^j(X_0,\calO(L_0)),
                \quad j=0,1,\ldots
\end{equation}
In \refss{prmainC}, we will prove that the restriction of $\p_t$
to $E^{*,K}_{\lam,t}$ vanishes for $t\gg0$. In other words, the
cohomology of the complex $(E^{*,K}_{\lam,t},\p_t)$ is isomorphic
to $E^{*,K}_{\lam,t}$. Hence, the inequalities \refe{Tzin} are, in
fact, equalities. This will complete the proof of \reft{mainC}.

Since we need some of the results obtained by Tian and Zhang
during their proof of \reft{TZ1}, we now recall briefly the main
steps of this proof.

\ssec{pr}{Main steps of the proof of \reft{TZ1}}
The proof of \reft{TZ1} in \cite{TianZh98} goes approximately as
follows: one shows first, that the eigenforms of $\lap^{*,K}$ with
eigenvalues smaller than $\lam$ are concentrated near
$\mu^{-1}(0)$. Then, using the local form of $\lap^{*,K}$ near
$\mu^{-1}(0)$, one describes the ``asymptotic behaviour" of
$E_{\lam,t}^{*,K}$ as $t\to\infty$.

More precisely, recall that $\oq:\mu^{-1}(0)\to X_0=\mu^{-1}(0)/K$
denote the quotient map and set
\[
      \tilh(x) \ = \ \sqrt{\Vol\, \oq^{-1}(x)}.
\]
Let $g^{L_0}$ and $g^{X_0}$ denote the Hermitian metric on $L_0$
and the Riemannian metric on $X_0$ induced by the fixed metrics on
$L$ and $X$ respectively (cf. \refss{lapt}). Set
$g^{L_0}_{\tilh}=\tilh^2g^{L_0}$ and let  $\p^*_\tilh$ denote the
formal adjoint of the Dolbeault differential $\p:\Ome^{0,*}(X_0,
L_0)\to \Ome^{0,*+1}(X_0, L_0)$ with respect to the metrics
$g^{L_0}_{\tilh}, g^{X_0}$. Let
\[
        \h{j} \ = \ \Ker\big(\,
                \p\p^*_\tilh+\p^*_\tilh\p:\, \Ome^{0,j}(X_0, L_0)\to\Ome^{0,j}(X_0, L_0)
                \big)
\]
be the space of harmonic forms. Then, for any $t\gg0$, one
constructs an isomorphism of vector spaces
\begin{equation}\label{E:map1}
        \Phi_{\lam,t}^j:\, \calH^j(X_0,L_0) \ \to \ E_{\lam,t}^{j,K},
                \qquad j=0,1,\ldots
\end{equation}
This implies the equality \refe{TZ1}.

Since we will use the above isomorphism in our proof of
\reft{mainC}, we now review briefly its construction and main
properties.

\rem{d_Q}
In \cite{TianZh98}, Tian and Zhang considered the operator
\[
        D_Q \ = \ \sqrt{2}\, (\,
                 \tilh\p\tilh^{-1}+\tilh^{-1}\p^*\tilh\, ):\,
                 \Ome^{0,j}(X_0,L_0) \ \to \ \Ome^{0,j}(X_0,L_0),
\]
where $\p^*$ denotes the formal adjoint of $\p$ with respect to
the metrics $g^{L_0},g^{X_0}$.  Then Tian and Zhang used the
method of \cite{BisLeb91} to construct a map from $\Ker D_Q$ to
$E_{\lam,t}^{j,K}$. Clearly, $\Ker
D_Q=\tilh^{-1}\calH^*(X_0,L_0)$. Our map \refe{map1} is a
composition of multiplication by $\tilh^{-1}$ with the map
constructed in \cite{TianZh98}.
\erem

\ssec{map1}{A bijection from $\h{*}$ onto $E_{\lam,t}^{*,K}$}
As a first step in the construction of the isomorphism
\refe{map1}, we construct an auxiliary map
$\Psi_{t}^j:\Ome^{0,j}(X_0,L_0)\to\Ome^{0,j}(X,L)$.

Let $N\to\mu^{-1}(0)$ denote the normal bundle to $\mu^{-1}(0)$ in
$X$. If $x\in\mu^{-1}(0),\ Y\in N_x$, let $t\in\RR\to
y_t=\exp_x(tY)\in X$ be the geodesic in $X$ which is such that
$y_0=x, \ dy/dt|_{t=0}=Y$. For $0<\eps<+\infty$, set
\[
        B_\eps \ = \ \{\, Y\in N: \ |Y|<\eps\, \}.
\]
Since $X$ and $\mu^{-1}(0)$ are compact, there exists $\eps_0>0$
such that, for $0<\eps<\eps_0$, the map $(x,Y)\in N \to\exp_x(tY)$
is a diffeomorphism from $B_{\eps}$ to a tubular neighborhood
${\calU}_\eps$ of $\mu^{-1}(0)$ in $X$. From now on, we will
identify $B_\eps$ with $\calU_\eps$. Also, we will use the
notation $y=(x, Y)$ instead of $y=\exp_x(Y)$.

Set $r=\dim K$. Since $0\in \grk^*$ is a regular value of $\mu$,
$\mu^{-1}(0)$ is a non-degenerate critical submanifold of
$|\mu|^2$ in the sense of Bott.  Thus, there exists an equivariant
orthonormal base $f_1,\ldots, f_r$ of $N$ such that, for any
$Y=y_1f_1 +\cdots + y_{r} f_{r}$,
\eq{mu2}
        |\mu(x,Y)|^2 \ = \
                \sum_{i=1}^{r}\, a_iy_i^2 \ + \ O(|Y|^3),
\end{equation}
where each $a_i$ is a positive $K$-invariant function on
$\mu^{-1}(0)$.

Let $p:\calU_\eps\to \calU_\eps/K$ denote the projection and let
$h$ be the smooth positive function on $\calU_\eps$ defined by
\[
        h(u) \ := \ \sqrt{\Vol(p^{-1}(p(u)))}, \quad u\in \calU_\eps.
\]
Note that, for any $u\in\mu^{-1}(0)$ we have $h(u)=\tilh(p(u))$,
cf. \refss{pr}.

The following simple lemma plays an important role in
\refss{integerb}:
\lem{h=a}
If $x\in\mu^{-1}(0)$, then $h(x)=(a_1(x)\cdots a_r(x))^{1/4}$.
\elem
\prf
Let $v_i,\mu_i$ be as in \refss{dt}. Let $J$ and $\<\cdot,\cdot\>$
denote the complex structure and the \ka scalar product on $TX$.
Then
\[
        \<v_i,v_j\> \ = \ \ome(v_i,Jv_j) \ = \ d\mu_i(Jv_j).
\]
Each $\mu_i$ is a function on $\calU_\eps$ and, hence, may be
written as $\mu_i=\sum \alp_{ij}(x)y_j$.

For any $x\in\mu^{-1}(0)$, consider the $r\times r$-matrices:
\[
        A(x) \ = \ \big\{\, \alp_{ij}(x)\, \big\}, \quad
        V(x) \ = \ \big\{\, \<v_i(x),v_j(x)\>\, \big\}, \quad
        F(x) \ = \ \big\{\, \<f_i(x),Jv_j(x)\>\, \big\}
\]
Clearly,
\[
        A F \ = \ V, \quad A^2 \ = \ \diag(a_1\nek a_r),
        \quad (\det F)^2 \ = \ \det V \ = \ h^4.
\]
Hence, $h^2=\det V/\det F = \det A = \sqrt{a_1\cdots a_r}$.
\eprf

For any $t>0$, consider the function $\bet_t$ on $\calU_\eps$
defined by the formula
\[
        \bet_t(x,Y) \ := \
            \Big(\, \prod_{i=1}^{r} a_i\, \Big)^{1/4}\,
             \Big(\, \frac{t}{2\pi} \,  \Big)^{r/4}\,
                \exp (-\frac{t}2\sum_{i=1}^r a_i y_i^2),
\]
where $x\in\mu^{-1}(0), Y\in N_x\simeq\RR^r$. Clearly,
\eq{bet}
    \int_{\RR^r}\, |\bet_t(x,Y)|^2\, dy_1\dots dy_r = 1,
\end{equation}
for any $x\in\mu^{-1}(0)$.

Let $\sig:X\to [0,1]$ be a smooth function on $X$, which is
identically equal to 1 on $\calU_{\frac\eps2}$ and such that
$\supp\sig\subset\calU_\eps$. We can and we will consider the
product $\sig\frac{\bet_t}{h}$ as a function on $X$, supported on
$\calU_\eps$.

Recall, from \refss{Mquot}, that $q:\Xss\to X_0$ denotes the
projection. Set
\[
        \Psi^j_t \ := \ \sig\frac{\bet_t}{h}\circ q^*:\,
                \Ome^{0,j}(X_0,L_0)\to\Ome^{0,j}(X,L),
        \qquad j=0,1,\dots
\]
Here $q^*:\Ome^{0,j}(X_0,L_0)\to\Ome^{0,j}(\Xss,\Lss)$ denotes the
pull-back, and we view multiplication by the compactly supported
function $\sig\frac{\bet_t}{h}$ as a map from
$\Ome^{0,j}(\Xss,\Lss)$ to $\Ome^{0,j}(X,L)$.

It follows from \refe{bet}  and the definition of $h$ that, for
$t\gg 0$, the map $\Psi_t^*$ is closed to isometry, i.e.,
\eq{psi*}
    \lim_{t\to\infty}\, \|\Psi_t^*\Psi_t-\Id\|=0.
\end{equation}

\rem{TZ1}
Clearly $\Psi_t^*$ does not commute with differentials, due to the
presence of the cut-off function $\sig$ in the definition of
$\Psi_t^*$ and also because $|\mu|^2$ is only approximately equal
to $\sum{}a_iy_i^2$. However, for any
$\alp\in\Ome^{0,*}(X_0,L_0)$, the restriction of $\Psi_t^*\p\alp$
to $\calU_{\frac\eps2}$ is ``very close" to $\p_t\Psi_t^*\alp$.
Since, for large values of $t$, ``most of the norm" of
$\Psi_t^*\p\alp$ is concentrated in $\calU_{\frac\eps2}$, we see
that, for $t\gg0$, the map $\Psi_t^*$ ``almost commutes" with
differentials. More precisely,
$\lim_{t\to\infty}\|\Psi_t^*\p\alp-\p_t\Psi_t^*\alp\|=0$ for any
$\alp\in\Ome^{0,*}(X_0,L_0)$.
\erem

Let $\oE^j_t\subset \Ome^{0,j}(X, L)$ denote the image of $\h{j}$
under $\Psi^j_t$.

The following theorem, which combines Theorem~3.10 and
Corollary~3.6 of \cite{TianZh98} and Theorem~10.1 of
\cite{BisLeb91}, shows that the image of $\h{*}$ under $\Psi^*_t$
is ``asymptotically equal" to $E^{*,K}_{\lam,t}$.
\th{TZ2}
Let $P_{\lam,t}^j: \Ome^{0,j}(X, L)\to E_{\lam,t}^{j,K}$ be the
orthogonal projection and let
$\Id:\Ome^{0,j}(X,L)\to\Ome^{0,j}(X,L)$ be the identity operator.
Then, there exists $\lam>0$, such that
\[
        \lim_{t\to\infty}\,
         \|(\Id-P_{\lam,t}^j)|_{\oE^j_t}\|
         \ = \  0,
\]
where $\|(\Id-P_{\lam,t}^j)|_{\oE^j_t}\|$ denotes the norm of the
restriction of \/ $\Id-P_{\lam,t}^j$ to $\oE^j_t$.
\eth

Define the map
\[
        \Phi_{\lam,t}^j:\, \calH^j(X_0,L_0) \ \to \ E_{\lam,t}^{j,K},
\]
by the formula
$\Phi_{\lam,t}^j\overset{\text{def}}{=}P_{\lam,t}^j\circ\Psi_{t}^j$.
It follows from \reft{TZ2}, that, the map $\Phi_{\lam,t}^j$ is a
monomorphism. With a little more work, cf. \cite{TianZh98}, one
proves the following
\th{TZ3}
There exists $\lam,t_0>0$ such that, for every $j=0,1,\dots$ and
every $t>t_0$, the map $\Phi_{\lam,t}^j:\calH^j(X_0,L_0)\to
E_{\lam,t}^{j,K}$ is an isomorphism of vector spaces.
\eth
\noindent This implies, in particular, \reft{TZ1}.

From \refe{psi*} and \reft{TZ2} we also obtain the following
\cor{phi-psi}
\(\displaystyle \lim_{t\to\infty}\|\Phi_{\lam,t}^*-\Psi_t^*\|=0.\)
\ecor

\sec{integr}{The integration map.  Proof of \reft{mainC}}

In this section we construct a map
$I_t:\Ome^{0,j}(X,L)\to\Ome^{0,j}(X_0,L_0)$, whose restriction to
$E_{\lam,t}^{*,K}$ is ``almost equal" to $(\Phi^*_{\lam,t})^{-1}$
for $t\gg0$, and such that $I_t\p_t=\p I_t$, for any $t\ge0$. The
very existence of such a map implies (cf. \refc{integr}) that the
restriction of $\p_t$ onto the space $E_{\lam,t}^{*,K}$ vanishes.
\reft{mainC} follows then from \reft{TZ1}.

\ssec{integr}{The integration map}
Recall that $q:\Xss\to X_0=\Xss/G$ is a fiber bundle. Recall,
also, that we denote $r=\dim_\RR K=\dim_\CC G$. Hence, $\dim_\CC
q^{-1}(x)=r$, for any $x\in{}X_0$. The action of $G$ defines a
trivialization of $\Lss$ along the fibers of $q$. Using this
trivialization, we define a map
$I_t:\Ome^{0,j}(X,L)\to\Ome^{0,j}(X_0,L_0)$ by the formula
\begin{equation}\label{E:integr}
        I_t:\, \alp \ \mapsto \
            \Big(\frac{t}{2\pi}\Big)^{r/4}\, \int_{q^{-1}(x)}\,
                e^{-\frac{t|\mu|^2}2}\alp\wedge\ome^r,
\end{equation}
Though the integral is taken over a non-compact  manifold
$q^{-1}(x)$ it is well defined. Indeed, by \cite[\S4]{GuiSter82},
$q^{-1}(x)$ is the set of smooth points of a  complex analytic
submanifold of $X$. Hence, cf. \cite[\S0.2]{GrifHar94}, the
Liouville volume $\int_{q^{-1}(x)}\ome^r$ of $q^{-1}(x)$ is
finite. It follows that the integral in \refe{integr} converges.
Moreover, it follows from \refe{mu2} that there exists a constant
$C>0$ such that
\eq{normI}
    \|\, I_t\, \| \ \le \ C,
\end{equation}
for any $t>0$.

The following theorem describes the main properties of the
integration map $I_t$.
\th{integr}
a. \ $\p\circ I_t = I_t\circ\p_t$, for any $t\ge0$.

b. \ Let $i:\h{*}\to \Ome^{0,*}(X_0,L_0)$ denote the inclusion and
let $\|I_t\circ\Phi_{\lam,t}^*-i\|$ denote the norm of the
operator $I_t\circ\Phi_{\lam,t}^*-i:\h{*}\to\Ome^{0,*}(X_0,L_0)$.
Then
\[
        \lim_{t\to\infty}\,
                \|I_t\circ\Phi_{\lam,t}^*-i\| \ = \ 0.
\]
\eth
We postpone the proof of \reft{integr} to the next section. Now we
will show how it implies \reft{mainC} (and, hence, also
Theorems~\ref{T:mainA} and \ref{T:mainB}).
First, we establish the following simple, but important corollary
of \reft{integr}.
\cor{integr}
Recall that a positive number $t_0$ was defined in \reft{TZ3}.
Choose $t>t_0$ large enough, so that
$\|I_t\circ\Phi_{\lam,t}^*-i\|<1$, cf. \reft{integr}. Then,
$\p_t\gam=0$ for any $\gam\in E_{\lam,t}^{*,K}$.
\ecor
\prf
Let $t$ be as in the statement of the corollary and let
$\gam\in{}E_{\lam,t}^{*,K}$. Then, $\p_t\gam\in E_{\lam,t}^{*,K}$.
Hence, it follows from \reft{TZ3}, that there exists
$\alp\in\h{*}$ such that $\Phi_{\lam,t}^*\alp=\p_t\gam$.

By \reft{integr}.a, the vector $I_t\p_t\gam=\p I_t\gam\in
\Ome^{0,*}(X_0,L_0)$ is orthogonal to the subspace $\h{*}$. Hence,
\[
        \|\alp\| \ \le \
          \|I_t\p_t\gam-\alp\| \ = \
                \|I_t\Phi_{\lam,t}^*\alp-\alp\|.
\]
Since, $\|I_t\circ\Phi_{\lam,t}^*-i\|<1$, it follows that
$\alp=0$. Hence, $\p_t\gam=\Phi_{\lam,t}^*\alp=0$.
\eprf

\ssec{prmainC}{Proof of \reft{mainC}}
We have already mentioned in \refss{lapt}, that $E_{\lam,t}^{*,K}$
is a subcomplex of $(\Ome^{0,*}(X,L),\p_t)$, whose cohomology is
isomorphic to $H^*(X,\calO(L))^K$. Since, by \refc{integr}, the
differential of this complex is equal to 0, we obtain
\[
    \dim{}H^j(X,\calO(L))^K=\dim{}E_{\lam,t}^{j,K}.
\]
\reft{mainC} follows now from \reft{TZ1}. \hfill$\square$

\sec{printegr}{Proof of \reft{integr}}

\ssec{printegr1}{Proof of \reft{integr}.a}
The first part of \reft{integr} is an immediate consequence of the
following two lemmas:
\lem{mu=0}
\(\displaystyle
        \int_{q^{-1}(x)}\,  \mu_i\p\mu_i\wedge\alp\wedge\ome^r
             \ = \  0,
\)
for any $i,j=0,1,\ldots$ and any  $\alp\in\Ome^{0,j}(X,L)$.
\elem
\prf
For any $l=0,1,\ldots$, let
$\Pi_{0,l}:\Ome^{*,*}(X_0,L_0)\to\Ome^{0,l}(X_0,L_0)$ denote the
projection. Then
\eq{pmu=dmu}
        \int_{q^{-1}(x)}\,  \mu_i\p\mu_i\wedge\alp\wedge\ome^r
            \ = \
                \Pi_{0,j+1}\, \int_{q^{-1}(x)}\,  \mu_i d\mu_i\wedge\alp\wedge\ome^r.
\end{equation}
Using \refe{dmu^2}, we obtain
\begin{multline}\label{E:dmu=}
        \int_{q^{-1}(x)}\,  \mu_i d\mu_i\wedge\alp\wedge\ome^r
             \ = \
                \int_{q^{-1}(x)}\,   \mu_i \iot_{v_i}\ome\wedge\alp\wedge\ome^r \\
         \ = \
                \frac1{r+1}\, \int_{q^{-1}(x)}\,  \mu_i\,
                        \iot_{v_i}(\ome^{r+1}\wedge\alp)
             \ - \
               \frac1{r+1}\, \int_{q^{-1}(x)}\, \mu_i \ome^{r+1}\wedge\iot_{v_i}\alp.
\end{multline}
The first summand in \refe{dmu=} vanishes since the vector $v_i$
is tangent to $q^{-1}(x)$. The integrand in the second summand
belongs to $\Ome^{r+1,r+j}(X,L)$. It follows, that
\[
        \int_{q^{-1}(x)}\,  \mu_i d\mu_i\wedge\alp\wedge\ome^r
         \ \in \ \Ome^{1,j}(X_0, L_0).
\]
The lemma follows now from \refe{pmu=dmu}.
\eprf

\lem{di=id}
\(\displaystyle
        \int_{q^{-1}(x)}\,
           e^{-\frac{t|\mu|^2}2}\p\alp\wedge\ome^r \ = \
                \p\, \int_{q^{-1}(x)}\, e^{-\frac{t|\mu|^2}2}\alp\wedge\ome^r,
\) for any $j=0,1\ldots$ and any $\alp\in\Ome^{0,j}(X,L)$.
\elem
\prf
Since the complement of $\Xss$ in $X$ has real codimension $\ge2$,
there exists a sequence $\alp_k\in\Ome^{0,j}(X,L), \ k=1,2\ldots$
convergent to $\alp$ in the topology of the Sobolev space
$W^{1,1}$, and such that $\supp(\alp_k)\subset\Xss$. Hence, it is
enough to consider the case when support of $\alp$ is contained in
$\Xss$, which we will henceforth assume.

Let $\bet\in\Ome^{n-j,n-j-1}(X_0,L^*_0)$, where $n$ is the complex
dimension of $X$ and $L^*_0$ denotes the bundle dual to $L_0$.

By \refl{mu=0},
\begin{multline}\notag
    \int_{q^{-1}(x)}\,  e^{-\frac{t|\mu|^2}2}\p\alp\wedge\ome^r
    \ = \
    \int_{q^{-1}(x)}\,
        \p\big(e^{-\frac{t|\mu|^2}2}\alp\wedge\ome^r\big)\\
    \ + \
    t\sum_{i=1}^r\,
        \int_{q^{-1}(x)}\,
             e^{-\frac{t|\mu|^2}2}\mu_i\p\mu_i\wedge\alp\wedge\ome^r
    \ = \ \int_{q^{-1}(x)}\,
        \p\big(e^{-\frac{t|\mu|^2}2}\alp\wedge\ome^r\big).
\end{multline}
Hence,
\begin{multline}\notag
        \int_{X_0}\,
          \Big( \int_{q^{-1}(x)}\,
            e^{-\frac{t|\mu|^2}2}\p\alp\wedge\ome^r\Big)\wedge\bet\\
        \ = \
          \int_{\Xss}\,
              \p\big(e^{-\frac{t|\mu|^2}2}\alp\wedge\ome^r\big)\wedge q^*\bet
        \ = \
          \int_{X}\,
            \p\big(e^{-\frac{t|\mu|^2}2}\alp\wedge\ome^r\big)
                                                \wedge q^*\bet\\
        \ = \
            (-1)^j\, \int_{X}\,
                e^{-\frac{t|\mu|^2}2} \alp\wedge\ome^r\wedge q^*(\p \bet)
             \ = \
                (-1)^j\, \int_{\Xss}\,
                    e^{-\frac{t|\mu|^2}2} \alp\wedge\ome^r\wedge q^*(\p \bet) \\
         \ = \
                (-1)^j\, \int_{X_0}\,
                  \Big( \int_{q^{-1}(x)}\,
                     e^{-\frac{t|\mu|^2}2}\alp\wedge\ome^r\Big)\wedge \p\bet
         \ = \
            \int_{X_0}\,
                \p\, \Big( \int_{q^{-1}(x)}\,
                        e^{-\frac{t|\mu|^2}2}\alp\wedge\ome^r\Big)\wedge \bet.
\end{multline}
\eprf

This completes the proof of \reft{integr}.a. \hfill$\square$

\ssec{integerb}{Proof of \reft{integr}.b}
Fix $\alp\in\h{*}$. Using \refe{normI} and \refc{phi-psi}, we
obtain
\[
    \|\, I_t\Phi^*_{\lam,t}\alp - I_t\Psi^*_{t}\alp\, \| \ \le \
        C\, \|\, (\Phi^*_{\lam,t}- \Psi^*_{t})\alp\, \| = o(1),
\]
where $o(1)$ denotes a form, whose norm tends to 0, as
$t\to\infty$. Hence,
\begin{multline}\notag
        I_t\Phi^*_{\lam,t}\alp \ = \ I_t\Psi^*_{t}\alp \ + \ o(1)
       \ = \
        \Big(\frac{t}{2\pi}\Big)^{r/4}\int_{q^{-1}(x)}\,
            e^{-\frac{t|\mu|^2}2}\sig\frac{\bet_t}{h} q^*\alp\wedge\ome^r
                                                      \ + \ o(1)\\
     \ = \
        \Big(\,
           \Big(\frac{t}{2\pi}\Big)^{r/4}\int_{q^{-1}(x)}\,
                e^{-\frac{t|\mu|^2}2}\sig\frac{\bet_t}{h}\, \ome^r\,
             \Big)\cdot\alp \ + \ o(1).
\end{multline}

Recall that in \refss{map1}, we introduced coordinates in a
neighborhood $\calU_\eps$ of $\mu^{-1}(0)$. Using this coordinates
and the definition of the function $h$, we can write
\[
     \int_{q^{-1}(x)}\,
        e^{-\frac{t|\mu|^2}2}\sig\frac{\bet_t}{h}\ome^r \ = \
            \int_{\RR^r}\,  e^{-\frac{t|\mu|^2}2}h\sig\bet_t dy_1\cdots dy_r.
\]
Hence, from \refe{mu2} and \refl{h=a}, we obtain
\begin{multline}\label{E:1+o}
        \Big(\frac{t}{2\pi}\Big)^{r/4}\int_{q^{-1}(x)}\,
            e^{-\frac{t|\mu|^2}2}\sig\frac{\bet_t}{h}\ome^r\\
        \ = \
             \Big(\, \frac{t}{2\pi} \, \Big)^{r/2}
             \Big(\, \prod_{i=1}^{r} a_i\, \Big)^{1/2}\,
        \int_{\RR^r}\, \sig
            \exp (-\frac{t}2\sum_{i=1}^r a_i y_i^2)\, dy_1\cdots
            dy_r + o(1)\\
            \ = \
               1 \ + \ o(1).
\end{multline}
Thus, $I_t\Phi^*_{\lam,t}\alp=\alp+o(1)$. In other words, the
operator $I_t\Phi^*_{\lam,t}$ converges to $i$ as $t\to\infty$ in
the strong operator topology. Since the dimension of $\h{*}$ is
finite, it also converges in the norm topology. \hfill$\square$

\providecommand{\bysame}{\leavevmode\hbox
to3em{\hrulefill}\thinspace}

\end{document}